\documentclass[review]{elsarticle}

\usepackage{lineno,hyperref, xcolor}
\usepackage[margin=1in]{geometry}

\makeatletter
\def\ps@pprintTitle{%
 \let\@oddhead\@empty
 \let\@evenhead\@empty
 \def\@oddfoot{}%
 \let\@evenfoot\@oddfoot}
\makeatother

\usepackage{amsmath}
\usepackage{float}
\usepackage{graphicx}
\usepackage{subfig}
\allowdisplaybreaks

\usepackage{algorithm}
\usepackage{algpseudocode}

\algtext*{EndFor}
\algtext*{EndIf}
\algtext*{EndFunction}
\algtext*{EndWhile}

\algrenewcommand\algorithmicforall{\textbf{foreach}}
\algrenewcommand\algorithmicindent{.8em}

\biboptions{authoryear}

\begin{document}

\begin{frontmatter}
\title{Amazon Last-Mile Delivery Trajectory Prediction Using Hierarchical TSP with Customized Cost Matrix}

\author{Xiaotong Guo, Baichuan Mo, Qingyi Wang \footnote{In alphabetical order by the last name}}
\address{Department of Civil and Environmental Engineering, Massachusetts Institute of Technology}
\address{77 Massachusetts Ave, Cambridge, MA, USA}

\begin{abstract}

In response to the Amazon Last-Mile Routing Challenge, Team Permission Denied proposes a hierarchical Travelling Salesman Problem (TSP) optimization with a customized cost matrix. The higher level TSP solves for the zone sequence while the lower level TSP solves the intra-zonal stop sequence. The cost matrix is modified to account for routing patterns beyond the shortest travel time. Lastly, some post-processing is done to edit the sequence to match commonly observed routing patterns, such as when travel times are similar, drivers usually start with stops with more packages than those with fewer packages. The model is tested on 1223 routes that are randomly selected out of the training set and the score is $0.0381$. On the 13 routes in the given model apply set, the score was $0.0375$.

\end{abstract}

\end{frontmatter}

\section{Introduction}

This report presents the thought processes, selected methodology, and expected results of the Amazon Last-Mile Routing Research Challenge by Team Permission Denied. In summary, the team went through four phases before arriving at the final submission. 

\textbf{Descriptive Analysis:} Upon receiving the challenge, a thorough descriptive analysis is done. The first important finding is that, in most circumstances, the drivers finish all deliveries in one zone before moving on to the stops in another zone. This rule is only broken when backtracking exists. A further look at the scores confirms this intuition: assuming the zone sequence and intra-zonal stop sequence are correct, the loss on the score due to certain zones being revisited is only 0.009. If the zone sequence is correct and the stops in each zone are shuffled, the average score is around 0.02. Therefore, getting the zone sequence correct is the most important, and the team decides to adopt a hierarchical approach: solving for the zone sequence, and then the intra-zonal stop sequence. This greatly reduces the scale of the problem since the majority of the routes have around 150 stops (up to 250), but the number of zones is between 6 and 47. Second, the zonal transitional probabilities are investigated. As most of the zones only appear in the training set once, an attempt at a frequency tabulation is not successful. On the other hand, 74\% of the zonal transitions select the zone that is closest by travel time, making the step-by-step prediction algorithm potentially successful. Next, the correlation between package dimensions, package counts, delivery time windows, and sequence order is investigated but no apparent relationship is found. 

\textbf{Benchmarking:} A benchmark model is created to establish an idea of the solution quality and expected performance. Since most drivers follow the navigation given by Amazon, a shortest-distance tour becomes a natural benchmark. The team solves a tour-based (where the start and end stations are both INIT) to generate zone sequences and a path-based (where the distance from the last zone to INIT is not counted) Travelling Salesman Problem (TSP) to generate intra-zonal stop sequences as benchmarks. Inside each zone, a path-based TSP is generated from the stop closest to the last zone to the stop closest to the next zone.

\textbf{Model Attempts:} Both naive TSP solutions achieve scores reasonable scores (around 0.06). To improve the performance, machine learning models are attempted. First, it is noticed that correctly predicting the first zone would significantly improve the TSP performance, therefore a neural network is constructed to predict the first zone based on the travel time, distance, package count and size, etc. Second, pure machine learning models to generate sequences are investigated, including myopic approaches that predict the next element based on previously predicted stops, as well as sequence-to-sequence (seq2seq) approaches that encode and decode the entire sequence. Third, different training methods are considered, including the traditional cross-entropy loss, customized weighted loss, as well as reinforcement learning using policy gradients. Lastly, some improvements are made to the benchmark TSP models by adding penalty costs to non-consecutive zone-ids. Due to the small sample size (6k), machine learning techniques cannot outperform the benchmark models. After experimenting with various modeling techniques, the team decides to use the TSP solution as the final submission.

\textbf{Hyperparameter Searching and Post-Processing:} The customized cost matrix involves hyperparameters that the team searched for over the given training set. Lastly, some post-processing patterns are identified to further improve the quality of our solution.

The highlights of the final submitted model are:
\begin{itemize}
    \item Hierarchical modeling - To reduce the size of each optimization problem, the problem is broken down into zone routing and intra-zonal stop routing.
    \item Customized TSP cost matrix - To account for considerations in addition to shortest distance, the cost matrix is modified and the TSP performance improved by almost 0.01.
    \item Post-processing to match behavioral patterns - Some TSP sequences are reversed to accommodate delivery patterns such as stops with more packages are visited first instead of last, all else being equal.
    \item Stable hyperparameters - The cost hyperparameters have good generalizability and do not require re-training.
\end{itemize}

The rest of the technical report reviews the relevant literature and its compatibility with the research question; describes the selected model in detail, and discusses the expected results.

\section{Literature Review}
This problem is essentially a vehicle routing problem, except that the traditional setup for vehicle routing problems aims for the shortest distance traveled, but the problem of interest looks for the most similarity with the observed sequence. Two research communities have extensively studied the vehicle routing problem: machine learning and operations research. Literature in both communities is reviewed, with the pros and cons of the algorithms discussed for the problem of interest.

\subsection{Operations Research}

Given a set of locations one would like to visit, a Traveling Salesman Problem (TSP) can be solved to find the route with the minimum cost or distance. The overview and history of the TSP can be found in \citet{TSP_overview}. 
Although TSP is a well-known NP-hard problem in combinatorial optimization, off-the-shelf integer optimization solvers (e.g., Gurobi and GLPK) are able to solve it efficiently for real-world instances.
One key approach we utilized when solving the TSP is the cutting-plane method~\citep{Marchand_Martin_Weismantel_Wolsey_2002}, which is initially applied to TSP by \citet{10.2307/166695}.

\subsection{Machine Learning}

Two types of architectures can be used to re-order the input sequence: step-by-step or sequence-to-sequence (seq2seq). Step-by-step prediction involves predicting the stops one by one, given the information from previous stops, as well as candidate stops. Since the information from candidate stops are crucial, feed-forward neural networks are not a good candidate since it does not attribute features to candidates. Instead, a feed-forward neural network with alternative-specific utility is adopted \citep{WANG2020234}. This architecture draws the connection between discrete choice models with neural networks and uses neural networks to generate the utility for each candidate, and the candidate with the highest 'utility' is chosen. A sequence is then formed by repeatedly feeding the selected stop into the algorithm to get the next stop until the end of the sequence is reached. The advantage of this algorithm is that it is at the stop level instead of the sequence level. Therefore, the sample size, which is critical for the success of machine learning algorithms, is significantly larger than the seq2seq models. The disadvantage of this algorithm is that it is myopic and only sees the next step candidates while making a selection.

In recent years, a lot of seq2seq prediction algorithms have been developed, mainly for natural language processing (NLP) tasks. Compared to step-by-step prediction, seq2seq models comprise an encoder and a decoder. All elements in the sequence are encoded before decoding starts, therefore a global view is attained. The architecture of encoder and decoder often involves variants of the recurrent neural networks (ex. long-short term memory networks) \citep{Sutskever2014}, or attention \citep{Vaswani2017}. Most seq2seq problems are considered with mapping one sequence to another, whereas the problem of interest is concerned with re-ordering the input sequence. Pointer network is proposed to solve this type of problem, where the decoder uses self-attention to point to one of the input elements \citep{Vinyals2015}. The authors used a pointer network to solve TSP and achieved similar performance to TSP solvers. One drawback of the original pointer network is that it is sensitive to the order of inputs. The authors, therefore, added another encoding module to eliminate this influence \citep{Vinyals2016}. However, in our experiments, this dependency can be leveraged by arranging the input set in a meaningful sequence to improve performance. For example, ordering the input stops according to the TSP sequence would accelerate model convergence and improve the score. However, in the papers presented above, 1M training samples were fed into the network. Given that the training set only contains 6000 routes, score improvements on TSP solutions are unsuccessful.

The original pointer network uses cross-entropy loss (supervised learning). In this problem, the cross-entropy loss is very inefficient due to the way the score is calculated, since the loss only considers the probability of the correct position, and the loss for predicting all other positions is the same. But the scoring function considers similarity in addition to correctness. The scoring function is not differentiable and cannot be directly used as the loss function and use gradient descent. An alternative training method is reinforcement learning based on policy gradients \citep{Ma2019, Bello2019}. Using the well-known REINFORCE algorithm, we can directly optimize the non-differentiable score function. Researchers have found that this method has the same sample efficiency and better generalizability for TSP problems compared to supervised learning \citep{Joshi2019}. However, training with reinforcement learning in this particular problem with the sample size and given information also does not outperform TSP solutions.

\subsection{Proposed Method}
Our proposed method is built upon the traditional TSP with a customized distance matrix that implicitly contains drivers' routing behaviors for the Amazon last-mile delivery. 
Compared to the existing TSP framework, which minimizes the total vehicle travel distance, we modified the distance matrix and generated optimal routes which minimized the total adjusted travel distance.

\section{Methodology}

\subsection{Data}\label{sec_data}

We observe that most of the drivers tend to visit all stops in a zone before going to the next zone. Hence, we divide the problem into two sub-problems. The first is to identify the zone sequence, and the second is to recognize the intra-zonal stop sequence. 

\begin{figure}[htb]
    \centering
    \includegraphics[width=1\linewidth]{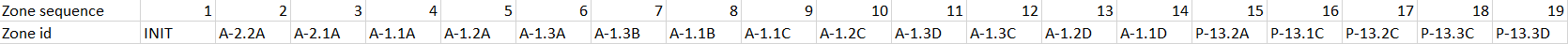}
    \caption{Example of zone sequence. ``INIT'' indicates the delivery station}
    \label{fig_zone_seq}
\end{figure}

The actual zone sequence is generated based on the order of each zone's first appearance. An example is shown in Figure \ref{fig_zone_seq}. For stops without zone id (due to missing data), we fill them with the zone ID of its (travel time-based) nearest stop. 

Three important properties are noticed while observing the zone sequences:
\begin{itemize}
    \item Most likely, the driver would finish a ``major zone'' first, then move to the next ``major zone''. A major zone is defined as the zone ID before the dot. For example, the major zone for ``A-2.2A'' is ``A-2''. For example, in Figure \ref{fig_zone_seq}, the driver first finishes major zone ``A-2'', then ``A-1'', finally ``P-13''. 
    \item Within a specific major zone, two adjacent ``inner zone'' ids are most likely have a ``difference of one''. The ``inner zone'' is defined as the zone ID after the dot. For example, the inner zone for ``A-2.2A'' is ``2A''. The ``difference of one'' is defined as follows. Given two inner zone IDs ``XY'' and ``AB'', where X and A are numbers and Y and B are characters, we have 
    \begin{align}
        |X - A| + |\texttt{ord}(Y) - \texttt{ord}(B)| = 1
    \end{align}
    where $\texttt{ord}(\cdot)$ function returns an integer representing the Unicode character. For example, ``1A'' and ``1B'' has a difference of one, so as ``1A'' and ``2A''. But ``1A'' and ``2B'' has a difference of two.
    \item When a driver finishes a ``major zone'' and move to another, the two adjacent major zone IDs are most likely to have a ``difference of one''. For example, in Figure \ref{fig_zone_seq}, the driver first finishes major zone ``A-2'', then ``A-1''. Those two major zone IDs have a difference of one.
\end{itemize}

To validate these three properties, we calculate the frequency that these rules hold in the data set. For all adjacent zone ID pairs, 87.67\% of them have the same major zone ID (Property 1). For all adjacent zone ID pairs within a specific major zone, 82.49\% of them have a ``difference of one'' (Property 2). For all adjacent zone ID pairs with major zone ID changes, 96.17\% of these changes lead to a ``difference of one'' between two major zone IDs (Property 3). These statistics support the three properties, which implies that \textbf{the zone ID includes a lot of information for the sequence estimation}. 

Another information we use is the planned service time and package volumes. Details on how these are utilized are shown in Section \ref{sec_post_process}.

We also collected outside data sources from OpenStreetMap. Specifically, we extract the number of traffic signals and highway ramps around every stop. Unfortunately, this does not help to improve our model, thus is dropped from our final submission. 

For the model's validation, we randomly separate the 6,112 routes into a training data set (4,889 routes) and a testing data set (1,223 routes), though our proposed solution does not require a training process.   




\subsection{Travelling Salesman Problem Formulation}

With the observation that drivers visit all stops within the same zone first and then move to the next zone, we solve a standard TSP with a modified travel time matrix to generate zone sequence first and then solve multiple path-TSP to identify intra-zonal stop sequence. 

First, we provide the formulation of the standard TSP solved for generating zone sequences. 
For a route instance with $n$ zones, the set of zones is indexed by $[n] = \{1,...,n\}$ and the initial station location is indicated by index $0$. 
Let $V$ represent the set of all locations that need to be visited including the initial station, i.e., $V = \{0, 1, ..., n\}$. 
$t_{ij}$ denotes the travel time between any two locations, i.e., $\forall i \neq j \in V$.
The travel time between any two zones is calculated as the average travel time between all possible pairs of stops between two zones. 
The decision variable for this problem is $x_{ij} \in \{0, 1\}, \; \forall i,j \in V$. $x_{ij} = 1$ indicates that the driver will visit to the location $j$ after visiting $i$.
Then, the TSP problem can be formulated as:

\begin{subequations}
\label{eq:Tour_TSP}
\begin{align}
    \min \quad & \sum_{i=0}^n \sum_{j=0}^n t_{ij} x_{ij} \\
    \text{s.t.} \quad
    & \sum_{i=0}^n x_{ij} = 1 \quad \forall j \in V \\
    & \sum_{j=0}^n x_{ij} = 1 \quad \forall i \in V \\
    & \sum_{i \in S} \sum_{j \notin S} x_{ij} \geq 1 \quad \forall S \subset V, S \neq \emptyset, V \\
    & \sum_{i \notin S} \sum_{j \in S} x_{ij} \geq 1 \quad \forall S \subset V, S \neq \emptyset, V \\
    & x_{ii} = 0 \quad \forall i \in V \\
    & x_{ij} \in \{0, 1\} \quad \forall i,j \in V
\end{align}
\end{subequations}

Where the objective (\ref{eq:Tour_TSP}a) minimizes the total travel time for the tour. 
Constraints (\ref{eq:Tour_TSP}b) and (\ref{eq:Tour_TSP}c) make sure that each visited location has exactly one predecessor and one successor in the optimal tour.
Constraints (\ref{eq:Tour_TSP}d) and (\ref{eq:Tour_TSP}e) are proposed to eliminate subtours in the optimal tour. 
Constraints (\ref{eq:Tour_TSP}f) avoid self loops and constraints (\ref{eq:Tour_TSP}g) guarantee decision variables are binary. 

The problem $(\ref{eq:Tour_TSP})$ is an Integer Linear Programming (ILP) with exponential number of constraints due to constraints (\ref{eq:Tour_TSP}d) and (\ref{eq:Tour_TSP}e). To solve this problem efficiently, we implemented both constraints (\ref{eq:Tour_TSP}d) and (\ref{eq:Tour_TSP}e) as lazy constraints, indicating they are only added to the problem if subtours are identified in the current optimal solution. 

To account for the observations made in the zone sequence (Section \ref{sec_data}), we propose three heuristics to modify the travel time matrix, which is the input for generating the optimal zone sequence.

\begin{enumerate}
    \item For travel time from the initial station to a zone $i$, if the zone is not within either i) $h$ closest zones from the initial station regarding travel times or ii) $h$ closest zones from the initial station regarding Euclidean distances, we modify the travel time to $t_{0i} * \alpha$, where $\alpha$ and $h$ are both parameters for the first proposed heuristic approach.
    \item For travel time between any two zones $i$ and $j$, if zone $i$ and zone $j$ are not from the same "major zone", we modify the travel time to $t_{ij} * \beta$, where $\beta$ is the parameter for the second proposed heuristic approach.
    \item For travel time between any two zones $i$ and $j$, if they are from the identical "major zone" and the difference between their zone ID after the dot does not equal to 1, we modify the travel time to $t_{ij} * \gamma$,  where $\gamma$ is the parameter for the third proposed heuristic approach.
\end{enumerate}

In the final submitted algorithm, we used the grid search approach to finalize values for all four heuristic parameters: $h = 9$, $\alpha = 1.04$, $\beta = 3.8$, $\gamma = 2.5$.

Solving the problem (\ref{eq:Tour_TSP}) with the modified travel time matrix leads to the optimal zone sequence\footnote{Without loss of generality, we can assume the sequence starts from the initial station indexed by $0$.} $S^* = (0, s_1,...,s_n)$, where $s_i$ indicates the $i$-th zone visited in the optimal sequence after departing from the initial station.
Then we solve the intra-zonal stop sequence using path-based TSP.
Given a set of locations $V$ need to be visited and the starting location $v_o$ and the ending location $v_d$, we can formulate the path-TSP problem as follows:

\begin{subequations}
\label{eq:Path_TSP}
\begin{align}
    \min \quad & \sum_{i=0}^n \sum_{j=0}^n t_{ij} x_{ij} \\
    \text{s.t.} \quad
    & \sum_{i=0}^n x_{ij} = 1 \quad \forall j \in V \setminus \{v_o, v_d\} \\
    & \sum_{j=0}^n x_{ij} = 1 \quad \forall i \in V \setminus \{v_o, v_d\}\\
    & \sum_{j \in V} x_{v_o j} = \sum_{i \in V} x_{i v_d} = 1 \\
    & \sum_{j \in V} x_{v_d j} = \sum_{i \in V} x_{i v_o} = 0 \\
    & \sum_{i \in S} \sum_{j \notin S} x_{ij} \geq 1 \quad \forall S \subset V, S \neq \emptyset, V \\
    & \sum_{i \notin S} \sum_{j \in S} x_{ij} \geq 1 \quad \forall S \subset V, S \neq \emptyset, V \\
    & x_{ii} = 0 \quad \forall i \in V \\
    & x_{ij} \in \{0, 1\} \quad \forall i,j \in V
\end{align}
\end{subequations}

The path-TSP problem (\ref{eq:Path_TSP}) is similar to the standard TSP problem (\ref{eq:Tour_TSP}) except that there will be no predecessors for the starting location $v_o$ and no successors for the ending location $v_d$, indicating by constraints (\ref{eq:Path_TSP}d) and (\ref{eq:Path_TSP}e).
The complete sequence is generated according to Algorithm \ref{alg1} based on generated zone sequence, where a heuristic parameter $k = 3$ is utilized in the final implementation.

\begin{algorithm}[!h]
\caption{Complete sequence generation based on the generated zone sequence. \\ Input: optimal zone sequence $S^* = (0, s_1,...,s_n)$, heuristic parameter $k$.}
\label{alg1}
\begin{algorithmic}[1]
\Function{CompletePathGeneration}{$S^*$}
    \State $S^*_{complete} \gets \{0\}$  \Comment{Initialize the complete sequence with the initial station}
    \For {$s_i = s_1,...,s_n$}
        \State Find the previous visited zone $s_{i-1}$ and the next visited zone $s_{i+1}$
        \State Calculate the average travel time between any stop $v \in s_i$ to all stops in zone $s_{i-1}$ and zone $s_{i+1}$
        \State Find $k$ nearest stops in zone $s_i$ regarding to zone $s_{i-1}$ as the set $M$
        \State Find $k$ nearest stops in zone $s_i$ regarding to zone $s_{i+1}$ as the set $N$
        \State Solve $k^2$ path-TSP (\ref{eq:Path_TSP}) between any pair of stops in $M \times N$.
        \State Let the path $S^*_{i}$ with the minimum travel time as the optimal sequence of zone $i$
        \State Append the sequence $S^*_{i}$ to the complete sequence $S^*_{complete}$
    \EndFor
    \State \textbf{return} $S^*_{complete}$
\EndFunction
\end{algorithmic}
\end{algorithm}

It is worth mentioning that all TSP instances are solved with the open-source ILP solver GLPK implemented with programming language Julia \citep{Bezanson2017} and optimization package JuMP \citep{DunningHuchetteLubin2017}.
After generating the complete stop sequence $S^*_{complete}$, we enter the post-processing stage to further improve sequence performances.



\subsection{Post-Processing}\label{sec_post_process}
After solving the stop sequence by TSP, we observe that most of the high-score (i.e., low performance) routes are due to partially or fully reverse of the sequence (i.e., a sequence A-B-C-D is erroneously estimated as D-C-B-A). Hence, we propose a post-processing method to correct the erroneous estimation due to reversal. 

We observe two properties from the data set:
\begin{itemize}
    \item Most of the drivers tend to serve the business areas first. The potential reason may be that it also takes a longer time to deliver packages in a business building. Serving them first can make the total service time more controllable at the end of the journey. Hence, we expect that the planned service time at the first several stops is larger than that of the last several stops.
    \item Most of the drivers tend to deliver large-size packages first. This may be because carrying large-size packages in the vehicle is not fuel-efficient. 
\end{itemize}

Based on these properties, for every generated stop sequence by TSP, we check whether we need to reverse it. Given a generated route $i$, let $p^+_i$ (resp. $p^-_i$) be the average planned service time of the first (resp. last) $p\%$ stops in route $i$. We will reverse route $i$ if
\begin{align}
    \frac{p^-_i}{p^+_i} \geq \theta,
    \label{eq_planned_time}
\end{align}
where $p$ and $\theta$ are hyperparameters representing the proportion of stops and a threshold. We set $p=15$ and $\theta=1.22$  based on cross-validation on the test set. Eq. \ref{eq_planned_time} means that in a generated sequence if the planned service time for the last several stops is too large, we may have the reversal error and need to correct it by reverse the whole sequence. 

After process by Eq. \ref{eq_planned_time}, we fixed all sequences that are already reversed. For the remaining sequences, we further check whether they need to be reversed based on package volumes. Specifically,  given a generated route $i$, let $v^+_i$ (resp. $v^-_i$) be the total package columns (depth$\times$width$\times$height) of the first (resp. last) 15\% stops in route $i$. We will reverse route $i$ if
\begin{align}
    \frac{v^-_i}{v^+_i} \geq \eta,
    \label{eq_planned_time}
\end{align}
where $\eta = 3$ is used.

After post-processing, a sequence validity check is performed. Specifically, we check whether the first stop of the estimated sequence is the delivery station, and whether the estimated sequence has the same stop IDs as the actual one. If either of these two criteria does not hold, we return a sequence by simply sort the stops using zone IDs, which ensures that stops with the same zone IDs are close to each other.

\section{Results and Conclusions}
\subsection{Performance}

Although the submitted formulation does not require model training, we have separated the given training set into the training (4889) and test set (1223) for self-evaluation of the machine learning models. Therefore, all self-evaluation is done over the test set. To reduce the evaluation time, we implemented the scoring function using Cython. Compared to the evaluation code in Python provided by the challenge host team, our implementation evaluates the same set of routes by using only one-third of the computation time. 



\begin{figure}[H]
    \centering
    \includegraphics[width=.6\linewidth]{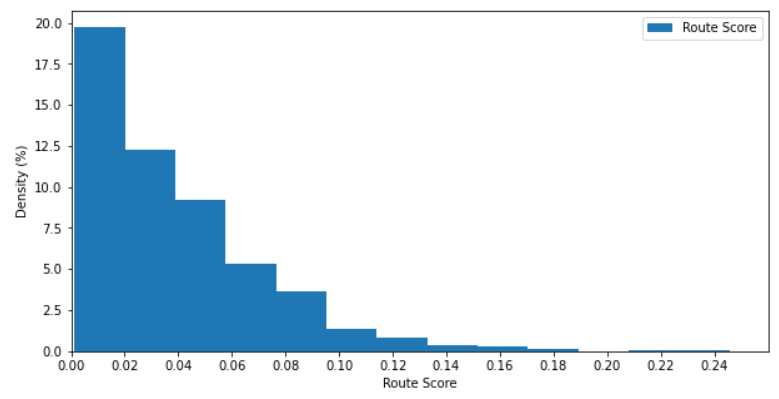}
    \caption{Route score performances.}
    \label{fig:route_score_performance}
\end{figure}

Figure \ref{fig:route_score_performance} shows the score distribution generated by our final algorithm. 
The route performance score follows an exponential distribution and most routes have a score below 0.1. 
The average route score is $0.0381$ for these 1223 testing routes. On the 13 routes in the given model apply set, the score was $0.0375$.

\subsection{Discussion}

\textbf{Zone sequence dominates score}. We observe that, if the zone sequence is perfectly predicted, even if the stop IDs within a zone are shuffled, the average route score can reach $0.0206$. Hence, most of our jobs focus on predicting the zone sequence, instead of the stop sequence. 

\textbf{The three properties of zone IDs (see Section \ref{sec_data}) may imply that drivers most likely follow the planned route and seldom deviate}. As the zone ID is used to ``help simplify the route planning process'' (quoted from Slack Q\&A), we believe that Amazon plans the route in a way that the zone IDs exhibit clear patterns. So the major challenge of this problem is to recover how Amazon plans the routes. This explains why TSP works better than machine learning methods under the given current information and sample size. 

\begin{figure}[H]
\captionsetup[subfigure]{justification=centering}
\centering
\subfloat[Actual route]{\includegraphics[width=0.5\linewidth]{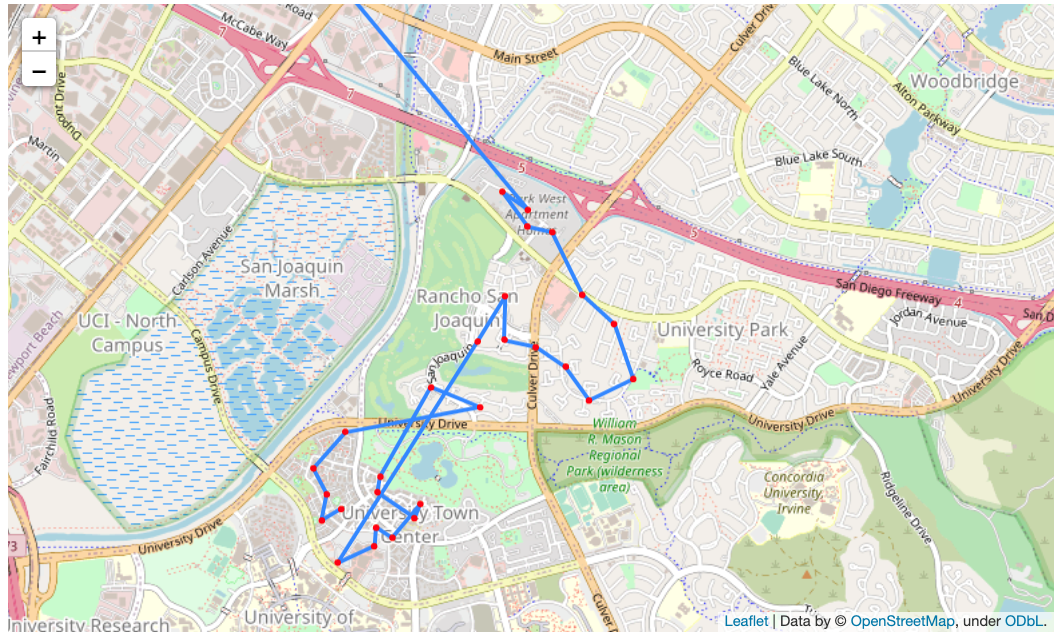}\label{sub_1}}
\subfloat[Predicted route by TSP]{\includegraphics[width=0.5\linewidth]{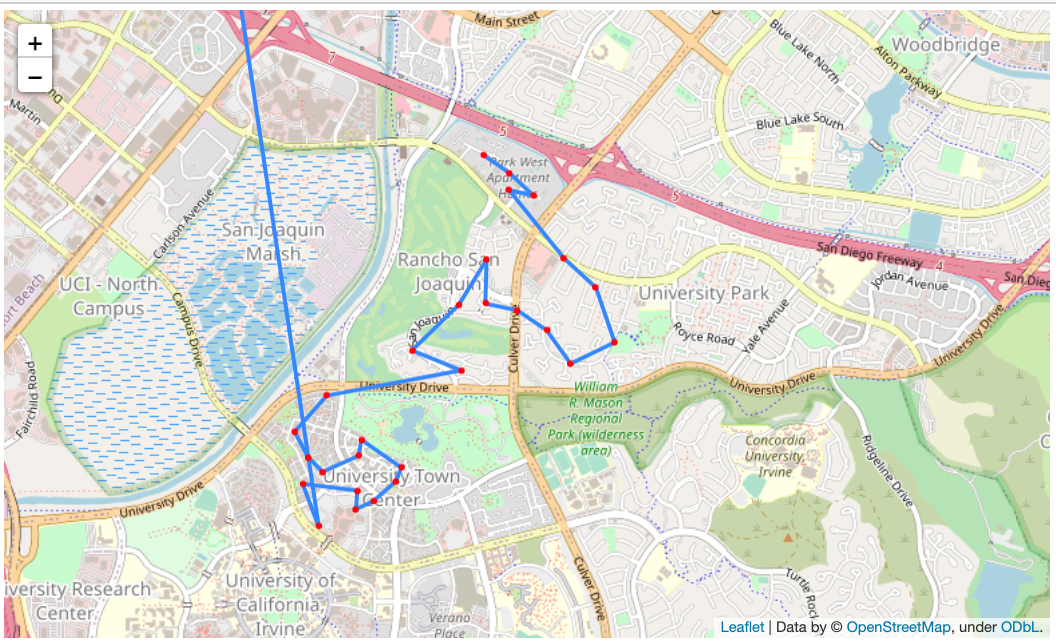}\label{sub_2}}
\caption{Examples of reverse prediction}
\label{fig_example_route}
\end{figure}

\textbf{The reversal problem remains}. Figure \ref{fig_example_route} shows an example of reverse prediction. Since we are not able to increase the first-zone prediction accuracy beyond 35\%, after post-processing, the reverse issues still exist. The post-processing reduces our score on our test set from 0.391 to 0.381. However, if we can have a 100\% correction rate for the reversal problems (i.e., always use the one with a smaller score), the score can reduce to 0.280, indicating that some further correction methods are needed. Note that we have tried to use the number of surrounding highway ramps as a new indicator, as well as using machine learning to predict the first zone, but it does not increase the model performance. 

\clearpage
\bibliography{mybibfile}

\end{document}